\documentclass[a4paper,11pt]{article}

\usepackage{amsmath}
\usepackage{amsthm}

\newtheorem{theorem}{Theorem}

\newtheorem{definition}[theorem]{Definition}

\newtheorem{question}{Question}

\DeclareSymbolFont{AMSb}{U}{msb}{m}{n}
\DeclareMathSymbol{\N}{\mathbin}{AMSb}{"4E}
\DeclareMathSymbol{\Z}{\mathbin}{AMSb}{"5A}
\DeclareMathSymbol{\R}{\mathbin}{AMSb}{"52}
\DeclareMathSymbol{\Q}{\mathbin}{AMSb}{"51}
\DeclareMathSymbol{\I}{\mathbin}{AMSb}{"49}
\DeclareMathSymbol{\C}{\mathbin}{AMSb}{"43}
\let\P\undefined
\DeclareMathSymbol{\P}{\mathbin}{AMSb}{`P}

\def\Belyui{Bely\u\i}

\begin{document}

\title{A filtration question on Bely\u{\i} pairs and dessins}

\author{Jonathan Fine\\Milton
  Keynes\\England\\\texttt{jfine@pytex.org}}

\date{28 September 2009}

\maketitle

\begin{abstract}
  A \Belyui{} pair is a holomorphic map from a Riemann surface to
  $S^2$ with additional properties.  A dessin d'enfants is a bipartite
  graph with additional structure.  It is well know that there is a
  bijection between \Belyui{} pairs and dessins d'enfants.

  Vassiliev has defined a filtration on formal sums of isotopy classes
  of knots.  Motivated by this, we define a filtration on formal sums
  of Bely\u{\i} pairs, and another on dessin d'enfants.  We ask if the
  two definitions give the same filtration.
\end{abstract}

\section{Introduction}

First, we recall some definitions~\cite{Belyui,Dessins}.  A
\emph{Bely\u{\i} pair} is a Riemann surface $C$ together with a
holomorphic map $f:C \to S^2 = \C \cup \{\infty\}$ to the Riemann
sphere, such that $f'(p)$ is non-zero provided $f(p)$ is not $0$, $1$
or $\infty$.  (Bely\u{\i} proved that given $C$ such an $f$ can be
found iff $C$ can be defined as an algebraic curve over the algebraic
numbers.)

A \emph{dessin d'enfants}, or \emph{dessin} for short, is a graph $G$
together with a cyclic order of the edges at each vertex, and also a
partition of the vertices $V$ into two sets $V_0$ and $V_1$ such that
every edge joins $V_0$ to $V_1$.  Necessarily, $G$ must be a bipartite
graph.  Traditionally, the vertices in $V_0$ and $V_1$ are coloured
black and white respectively.

It is easy to see that a Bely\u{\i} pair gives rise to a dessin, where
$V_0=f^{-1}(0)$, $V_1 = f^{-1}(1)$, and the edges are the components
of the inverse image $f^{-1}([0,1])$ of the unit interval in $\C$.
The cyclic order arise from local monodromy around the vertices.

A much harder result, upon which our definitions rely, is that up to
isomorphism every dessin arises from exactly one Bely\u{\i} pair, or
in other words that there is a bijection between isomorphism classes
of Bely\u{\i} pairs and dessins.

\section{Definitions}

\begin{definition}[Bely\u{\i} object]
  A \emph{Bely\u{\i} object} $B$ consists of $((B_C, B_f), B_D)$ where
  $(B_C, B_f)$ is a Bely\u{\i} pair and $B_D$ is the associated dessin
  (or vice versa for the dessin and the pair).
\end{definition}

\begin{definition}[Vassiliev space]
  The \emph{Vassiliev space} $V=V_\C$ (for Bely\u{\i} objects) is the
  vector space over $\C$ which has as basis the isomorphism classes of
  Bely\u{\i} objects.
\end{definition}

Clearly, when an edge is removed from a dessin then it is still a
dessin.  Suppose $D$ is a dessin, and $T$ is a subset of its edges.
We will use $D \setminus T$ to denote the dessin so obtained.  This
same operation can also be applied to a Bely\u{\i} object $B$, even
though computing the associated curve $(B\setminus T)_C$ from $B_D$
and $T$ might be hard.

We will now define one or two filtrations of $V$.

\begin{definition}[Dessin with $d$ optional edges]
  Let $D$ be dessin and $S$ a $d$-element subset of $D$. Each subset
  $T$ of $S$ determines a dessin $S\setminus T$ and hence a \Belyui{}
  object $B_{S\setminus T}$.  Let $|T|$ denote the number of edges in
  $T$.  Use
  \[
   B_S = \sum\nolimits _{T\subseteq S} (-1)^{|T|}B_{S \setminus T}
  \]
  to define a vector $B_S$ in $V$, which we call \emph{the expansion
    of a dessin with $d$ optional edges}.
\end{definition}

\begin{definition}[Dessin filtration]
   Let $V_{D,d}$ be the span of the expansions of all dessins with $d$
   optional edges.  The sequence
  \[
    V =
    V_{D, 0} \supseteq
    V_{D, 1} \supseteq
    V_{D, 2} \supseteq
    V_{D, 3} \ldots
  \]
  is the \emph{dessin filtration} of $V$.
\end{definition}

We can also think of a Bely\u{\i} object as a map $f:C\to S^2$ (with
special properties).  Let $(C_1, f_1)$ and $(C_2, f_2)$ be Bely\u{\i}
pairs.  Then there is of course a map
\[
  g: C_1 \times C_2 \to S^2 \times S^2 \>.
\]

Let $\Delta \subset S^2 \times S^2$ denote the diagonal, and let $C$
denote $g^{-1}(\Delta)$, and $f$ the restriction of $g$ to $C$.  In
general
\[
  f: C \to \Delta \cong S^2 
\]
will not be a Bely\u{\i} pair.  There are two possible problems.  The
first is that $C\subset C_1\times C_2$ might have self intersections
or be otherwise singular.  If this happens, we replace $C$ by its
resolution, which is unique.

The second problem is more interesting.  It might be that $f$ has
critical points not lying above the special points $0$, $1$ and
$\infty$.  This problem cannot be avoided.  However, the above
discussion does show that there is product, which we will denote by
`$\circ$', on holomorphic branched covers of $S^2$.

\begin{definition}[Product filtration]
  Let $W$ be the vector space with basis isomorphism classes of
  branched covers of $S^2$.  We set $W_n$ to be the span of all
  products of the form
  \[
  (A_1 - B_1) \circ
  (A_2 - B_2) \circ
  \ldots \circ
  (A_n - B_n)
  \]
for $A_i$ and $B_i$ basis vectors of $W$.  Clearly, the $W_n$ provide
a filtration of $W$.
\end{definition}

\begin{definition}[\Belyui{} filtration]

The induced filtration of $V$ defined by $V_{B,n} = W_n \cap V$ is
called the \emph{\Belyui{} filtration} of $V$.
\end{definition}

\section{Questions}

\begin{question}
Are the two filtrations $V_D$ and $V_B$ equal?
\end{question}

If so, then we have also answered the next two questions.

\begin{question}
  The absolute Galois group acts on \Belyui{} pairs, and preserves the
  \Belyui{} filtration.  Does this action also preserve the dessin
  filtration?
\end{question}

\begin{question}
  Because the dessins with $d$ edges, all of which are optional, span
  $V_d/V_{d+1}$, the dessin filtration has finite dimensional
  quotients.  Does the \Belyui{} filtration have finite dimensional
  quotients?
\end{question}

Investigating the last two questions might help us answer the first.
They might also be of interest in their own right.

\end{document}